\documentstyle[amssymb]{amsart}

\newcommand{\stopproof}{\hfill \nobreak\medskip $\blacksquare$ \\
\hspace*{\fill}}

\newcommand{\AND}{\mbox{ \rm and }}

\newcommand{\IF}{\mbox{ \rm if }}

\newcommand{\forces}[2]{\Vdash_{#1} \mbox{``} #2 \mbox{''}}

\newcommand{\proof}{{\bf Proof:} \ }
\newcommand{\labelx}[1]{\label{#1}\marginpar{#1}}

\newcommand{\lcard}{\, \mid \!}
\newcommand{\rcard}{\! \mid \,}

\newcommand{\card}[1]{\lcard #1 \rcard}
\newtheorem{theor}{Theorem}[section]

\newtheorem{defin}{Definition}[section]
\newtheorem{corol}{Corollary}[section]
\newtheorem{lemma}{Lemma}[section]

\newtheorem{claim}{Claim}
\author{Juris Stepr\={a}ns}
\address{Department of Mathematics\\ York University\\ Toronto,
Canada\ \ \ \ \ M3J 1P3 \\{\sc steprans\@nexus.yorku.ca}}

\newcommand{\Reals}{{\Bbb R}}
\newcommand{\Sacks}{{\Bbb S}}
\newcommand{\Model}{{\frak M}}
\newcommand{\Poset}{{\Bbb P}}
\newcommand{\cov}{\mbox{cov}}
\title{Decomposing with Differentiable Functions }
\thanks{The  author  is partially
supported by NSERC}

\begin{document}
\maketitle
\bibliographystyle{amsplain}
\begin{abstract}
Let ${\frak d}_{m,n}$ denote the least number of sets into which
$m$-dimensional Euclidean space can be decomposed so that each member
of the partion has an $n$-dimensional tangent plane at each of its
points. It will be shown to be consistent that ${\frak d}_{m,n}$ is
different from ${\frak d}_{m,n+1}$. These cardinal will be shown to be
closely related to the invariants associated with the problem of
decomposing continuous functions into differentiable ones.

\end{abstract}
\section{Introduction}
This paper is concerned with problems arising from partitioning sets
into subsets with differentiability properties. If $V$  is a vector
space and $v\in V$ then $\langle v\rangle$ will be used to dentoe the space
spanned by $v$.

\begin{defin}
Let $P_{k,n}$ denote the space of all subspaces of $\Reals^n$ of
dimension less than or equal to $k$. Define a metric $\rho$ on
$P_{k,n}$ by letting $$\rho(V,W) = \max_{x\in V^+}\min_{y\in W^+}
\arccos\left(\frac{x\cdot y}{2\|x\|\|y\|}\right)$$ where $V^+$ and $W^+$ are the
non-zero elements of the corresponding subspaces. If $v\in \Reals^n$
then $\rho(\langle v\rangle),V)$ will be abbreviated to $\rho(v,V)$
and if $v'\in
\Reals^n$ then $\rho(\langle v\rangle),\langle v'\rangle)$ will be
abbreviated to $\rho(v,v')$. 
\end{defin}

The following definition plays a central role in this paper.

\begin{defin}
If $A \subseteq \Reals^n$, $V \in P_{k,n}$ \labelx{TangentPlane}and $
b\in \Reals^n$ then 
$V$  will be said to 
be a $k$-dimensional tangent plane to $A$ at $b$ if and only if for
every $\epsilon > 0$ there is 
some $\delta > 0$ such that $\rho(V,a) < \epsilon$ for every $a
\in  A$ such that $0 < \|a-b\| < \delta$.
A subset $A \subseteq \Reals^n$ will be said to be a $k$-smooth
surface if for every $b \in A$ there is $V\in P_{k,n}$ which is a
$k$-dimensional tangent plane to $A$ at $b$.
\end{defin}

Observe that if $A\subseteq \Reals^n$ and $V$ is a $k$-dimensional
tangent plane to $A$ at $b$ and if there is no $m$-dimensional tangent
plane to $A$ at $b$ then $V$ is the unique $k$-dimensional tangent
plane to $A$ at $b$. However it is possible that a $k$-dimensional
smooth surface in $\Reals^n$ may have $m$-dimensional tangent planes
at some --- or even all --- points, where $m < k$.

\begin{lemma}
If $S\subseteq \Reals^m$\labelx{??} is an $n$-smooth surface then there is a
Borel set $S' \supseteq S$ which is also an $n$-smooth surface.
\end{lemma}
\proof Choose $D$, a countable dense subset $S$ and let $S'$ be the
set of all $x\in \Reals^m$ such that $D$ has an $n$-dimensional tangent
plane at $x$. \stopproof

\begin{defin}
Let ${\cal D}_{m,n}$ denote the $\sigma$-ideal generated by the
$n$-smooth surfaces of $\Reals^m$.
\end{defin}
The covering number of ${\cal D}_{m,n}$ will be denoted by ${\frak
d}_{m,n}$. Unanswered questions about the additivity, cofinality and
other invariants of ${\cal D}_{m,n}$ naturally suggest themselves.
\section{Lower Bounds}
For the sake of the next theorem, let ${\cal D}_{m,0}$ denote the
$\sigma$-ideal generated by the
discrete subsets of $\Reals^m$ and note that ${\frak
d}_{m,0} = 2^{\aleph_0}$.
\begin{theor} If $m > n\geq 1$ then
${\frak
d}_{m,n}^+ \geq {\frak
d}_{n,n-1}$. 
\end{theor}
\proof Let ${\frak d}_{m,n} =  \kappa$ and suppose that  $\kappa^+ < {\frak
d}_{n,n-1}$ . Let  $\{S_\alpha\}_{\alpha\in
\kappa}$ be $n$-smooth surfaces such that 
$\bigcup_{\alpha\in
\kappa}S_\alpha = \Reals^m$.
Since $m > n$ it is possible to let $V\subseteq \Reals^m$ be an
$n$-dimensional subspace and let $W$ 
be a 1-dimensional subspace orthogonal to $V$. Noting that $\kappa^+
< {\frak d}_{n,n-1}\leq 2^{\aleph_0}$, let $\{w_\xi\}_{\xi\in
\kappa^+}$ be distinct elements of $W$. For each $\xi \in \kappa^+$ let
$S_{\xi,\alpha}$ consist of all those $s\in S_\alpha\cap (V+w_\xi)$
such that $V$ does not contain  a tangent plane to $S_\alpha$ at $s$.
(Recall that tangent planes have been defined as subspaces so it does
not make sense to say that  $V+ w_\xi$ does not contain  a tangent
plane to $S_\alpha$ at $s$.) It follows that for each $s\in
S_{\xi,\alpha}$ the  plane tangent to $S_\alpha$ at $s$ intersects $V$
on a subspace of dimension smaller than $n$ and so 
$S_{\xi,\alpha}$,  considered as subspace of $V$ under the
obvious isomorphism, is a $(n-1)$-smooth surface. Since $\kappa^+ <
 {\frak d}_{n,n-1}$ it follows that there is some $v\in V$
such that $v+w_\xi\notin S_{\xi,\alpha}$ for all $\xi \in \kappa^+$ and
$\alpha \in \kappa$.  There is an uncountable
set $X\subseteq \kappa^+$ and $\alpha \in \kappa$ such that $v +
w_\xi\in S_\alpha$ for each $ \xi \in X$. Let $\xi \in X$ be such that
$w_\xi$  is a limit of $\{w_\zeta\}_{\zeta \in X}$. Since $v+w_\xi \in
S_\alpha \setminus S_{\xi,\alpha}$ it follows that $S_\alpha$ has a
tangent plane at $v+w_\xi$ which is contained in $V$. This contradicts
the fact that $W$ is orthogonal to $V$ because, since $v+w_\xi$ is a
limit of   $\{v+w_\zeta\}_{\zeta \in X}$, it must be that any tangent plane
at $v+ w_\xi$ includes $W$.
\stopproof

It is worth noting that if $k < n$ and $S\subseteq \Reals^n$ is
$k$-smooth then $S$ is nowhere dense. Hence 
the covering number of ${\cal D}_{m,n}$ is at least as great as that
of the meagre ideal.

\section{Partitioning the plane into a small number of differentiable sets}

Sacks forcing will be denoted by $\Sacks$. The following definition
generalizes Sacks forcing by introducing a forcing partial order
which, it will be shown, is intermediate between the product an
interation of a finite number of Sacks reals. Before giving the
definition, define  an indexed set $\{x_i\}_{i=0}^{n}\subseteq
\Reals^m$ to be $\epsilon$-orthogonal if and only if $\|\rho(x_i - x_0,
x_j- x_0) - \pi/2\| < \epsilon$ whenever $1 \leq i < j \leq
n$. Observe that elementary algebra implies that for each dimension $m$
there is $\epsilon(m) > 0$ such that if $\{x_i\}_{i=0}^{ n}\subseteq
\Reals^m$ is $\epsilon(m)$-orthogonal then 
for any subspace $V$ of dimension $k < n$ there is $i$ such that $1
\leq i \leq n$ and $\rho(x_i - x_0,V) > \epsilon(m)$.

\begin{defin}
The partial order \labelx{ModifSacks} 
$\Sacks(m,n)$ is defined to consist of all Borel
sets $B\subseteq \Reals^m$ such that there exists a family
$\{C_{\xi} : \xi \in {}^{\stackrel{\omega}{\smile}}(n+1)\}$ such that 
\begin{enumerate}
\item $C_\xi$ is an open ball of radius less than
$1/\card{\xi}$ 
\item $\bigcap_{i\in \omega}\bigcup_{\xi: i \to n+1}
\overline{C_{\xi}} \subseteq B$ 
\item $\overline{C_{\xi}}\cap \overline{C_{\xi'}} = \emptyset$ if
$\card{\xi} = 
\card{\xi'}$  and $\xi\neq\xi'$
\item $\overline{C_\xi}\subseteq C_{\xi'}$ if $\xi' \subseteq \xi$
\item if $x_k\in C_{\xi\wedge k}$ then $\{x_i\}_{i=0}^{ n}$ is $\epsilon$-orthogonal. 
\end{enumerate}
The family $\{C_{\xi} : \xi \in {}^{\stackrel{\omega}{\smile}}(n+1)\}$
will be called a witness to the fact that $B\in \Sacks(m,n)$. 
\end{defin}
If $B \in \Sacks(m,n)$ is such that there exist open balls $C_\xi$
such that $$\bigcap_{i\in \omega}\bigcup_{\xi: i \to n+1} C_{\xi} = B$$
and 
$\{C_{\xi} : \xi \in {}^{\stackrel{\omega}{\smile}}(n+1)\}$ is a
witness to the fact that $B\in\Sacks(m,n)$
 then $B$ will be referred to as a {\em natural}
member of $\Sacks(m,n)$. It is obvious that
the natural members of $\Sacks(m,n)$ are dense in $\Sacks(m,n)$.

\begin{lemma} If $k < n$,
$B$ is a\labelx{qqq} natural member of $ \Sacks(m,n)$ and $S$ is a
$k$-dimensional smooth surface in $\Reals^m$ then $S\cap B$ is meagre
relative to $B$.
\end{lemma}
\proof Since $B$ is a natural member of $ \Sacks(m,n)$,
there
exist open balls 
$\{C_{\xi} : \xi \in {}^{\stackrel{\omega}{\smile}}(n+1)\}$ such that 
$\bigcap_{i\in \omega}\bigcup_{\xi: i \to n+1} C_{\xi} = B$.
 For each $x\in
B$ let $\psi_x :\omega \to n+1$ be the unique function  satisfying
that $x \in C_{\psi_x\restriction j}$ for each $j\in \omega$. 

By Lemma \ref{??} it may as well be assumed that $S$ is Borel and,
hence, satisfies the Property if Baire. If $S$ is not meagre relative
to $B$ then let $U$ be an open set in $B$ such that
$B\cap U\cap S$ is comeagre in $B\cap U$. Notice that the set of all
$x \in B$ such that $\psi_x(j) = 0$ for only finitely many $j$ is a
meagre set in $B$. Let $x\in B\cap U\cap S$ be such that $\psi_x(j)
= 0$ for infinitely many $j \in \omega$. Let $V$ be a $k$-dimensional
tangent plane to $S$ at $x$ and let $\delta$ be such that
$\rho(V,a) < \epsilon(m)$ for every $a\in A$ such that $0 < \|a-x\|<\delta$.
   
Let $j\in\omega$ be such that $j > 1/\epsilon(m)$,  $j > 1/\delta$,
the open ball of radius $1/j$ around $x$ is contained in $U$ and
$\psi_x(j) = 0$. It follows that $C_{\psi_x\restriction 
j}\subseteq U$ and hence $C_{(\psi_x\restriction j)\wedge i} \subseteq
U$ for each $ i \in n$.  It follows from Definition~\ref{ModifSacks} that
$C_{(\psi_x\restriction j)\wedge i} \cap B$ is non-empty and open in
$B$. Hence it is possible to choose $x_i \in C_{\psi_x\restriction
j\wedge i} \cap B\cap S $ for each $ i \in n$ such that 
$\{x_i\}_{i=0}^n$ is an $\epsilon(m)$-orthogonal family.
Since the dimension of $V$ is less than $n$ it follows from the choice
of $\epsilon(m)$ such that $\rho(V,x_i - x_0) > \epsilon(m)$ for some
$i$ between 1 and $n$. This is a contradiction.
\stopproof    

\begin{lemma} If $B$\labelx{GDelta} is a natural member of
$\Sacks(m,n)$ and $W$ is 
a dense $G_\delta$ relative to $B$ then $W \in \Sacks(m,n)$.
\end{lemma}
\proof The proof is standard.
Let $C_\xi$ be open balls witnessing that $B$ is natural member of
$\Sacks(m,n)$.  
Let $W = \bigcap_{i\in\omega}U_i$ where each $U_i$ is a dense open
set.

Inductively choose $C'_\xi\subseteq U_{\card{\xi}}$ such that
$C'_\xi = C_\eta$ for some $\eta : k \to n+1$ where
$k \geq \card{\xi}$. Having chosen $C'_\xi$, use the fact that
$C'_\xi = C_\eta$ implies that $C'_\xi$ is open to choose
some  $U\subseteq C'_\xi\cap U_{\card{\xi}+1}$ which is open and
nonempty relative to $B$. This means that there is some $\mu \supseteq
\eta$ such that $C_\mu\cap B \subseteq U$. Let $C'_{\xi\wedge i} =
C_{\mu\wedge i}$ for $ i\in n+1$. It follows that the family of sets
$C'_\xi$ witnesses that $W\in \Sacks(m,n)$.
\stopproof

\begin{corol} If
$B\in \Sacks(m,n)$ and $S\in {\cal D}_{m,k}$ \labelx{TakeAway} then
$B\setminus S\in
\Sacks(m,n)$.\end{corol}
\proof Without loss of generality, assume that $B$ is natural. By
Lemma \ref{qqq} and the definition of ${\cal D}_{m,k}$, it follows
that $B\setminus S$ is comeagre relative to $B$. Now apply Lemma
\ref{GDelta}. \stopproof

\begin{corol}If
$B\in \Sacks(m,n)$ and\labelx{www} $B \subseteq \bigcup_{k\in\omega} A_{k}$
  where each $A_k$ is Borel then there is some
$k\in\omega$ such that $B\cap A_k \in \Sacks(m,n)$.
\end{corol}
\proof  Without loss of generality, assume that $B$ is natural. For
some $k$ it must be that $A_k$ is of second category relative to $B$.
Let $U$ be a non-empty open set such that $A_k$ is comeagre in
$B\cap U$. Then it is clear that $B\cap U \in \Sacks(m,n)$ so applying
Lemma \ref{GDelta} gives the desired result.
\stopproof

\begin{corol}
If $g\in \Reals^m$ is the \labelx{GoUp} generic real added by forcing with
$\Sacks(m,n)$ and $k < n$ then $g$ does not belong to any
$k$-smooth surface with a Borel code in the ground model.
\end{corol}
\proof Suppose that $S$ is a $k$-smooth surface and $B\in
\Sacks(m,n)$ is such that $$B\forces{\Sacks(m,n)}{g\in S}$$ From
Corollary \ref{TakeAway} it follows that $B\setminus S\in
\Sacks(m,n)$. However  $B\setminus S\forces{\Sacks(m,n)}{g\notin S}$.
\stopproof

\section{The iteration}
The countable
support iteration of length $\alpha$ of forcings 
$\Sacks(m,n)$ will be denoted by $\Poset_\alpha^{m,n}$. As a
convenience, by countable support will be meant that if $p \in  
\Poset_\alpha^{m,n}$ then $p(\gamma) = \Reals^m$ for all but countably
many $\gamma \in \alpha$.  From Corollary
\ref{GoUp} of the previous section it follows that if $V$ is a model
of the Generalized Continuum Hypothesis and
$G$ is  $\Poset_{\omega_2}^{m,n}$ generic over $V$ then ${\frak d}_{m,k}
= \omega_2$ for each $ k < n$ in $V[G]$. It will be shown in this
section that ${\frak d}_{m,n}
= \omega_1$ in $V[G]$. 

To this end, suppose that $V$ is a model
of the Continuum Hypothesis, $p^* \in \Poset_{\omega_2}^{m,n}$ and that 
$p^*\forces{\Poset_{\omega_2}^{m,n}}{x\in\Reals^m}$. It follows that
there must be some $\alpha \in \omega_2$ such that 
$p^*\forces{\Poset_{\alpha}^{m,n}}{x\in V[G\cap
\Poset_{\alpha}^{m,n}]}$ and
 $p^*\forces{\Poset_{\alpha}^{m,n}}{x\notin V[G\cap
\Poset_{\beta}^{m,n}]}$ for all $\beta \in \alpha$. It must be shown
that there is, in the ground model, an $n$-smooth surface $S
\subseteq\Reals^m$ as well as a condition $q
\leq p^*$ such that  $q\forces{\Poset_{\omega_2}^{m,n}}{x\in S}$.

 Fix $x$ and let $\Model \prec H(\aleph_3)$ be a countable elementary
submodel containing $p^*$, $x$ and $\alpha$. Notice that if $p$ and
$q$ belong to $\Poset_{\alpha}^{m,n}\cap \Model$ then the $p\wedge q $
also belongs to $\Model$. The exact definition of $p\wedge q$ is not
important here, only that it is definable in $\Model$. Observe that
$p\wedge q$ exists in $\Poset_{\alpha}^{m,n}$ because each factor of
the iteration has a naturally defined meet operation.

Let ${\cal P}(\Poset_{\alpha}^{m,n}\cap  \Model)$ be given the
Tychonoff product topology --- in other words, it is homeomorphic to the
Cantor set. Let ${\cal F}$ be the set of all $F\in {\cal
P}(\Poset_{\alpha}^{m,n}\cap  \Model)$ such that
\begin{itemize}
\item if $p \in F$ and $q \geq p$ then $ q\in F$
\item if $ p\in F$ and $q \in F$ then $p\wedge q \in F$ and $p\wedge q
\neq \emptyset$
\item if $p \in F$, $\{q_i\}_{i\in k}\subseteq \Poset_{\alpha}^{m,n}$
and there does not exist any $ r \in  \Poset_{\alpha}^{m,n}$ such that
$r \leq p$ and $r \wedge q_i = \emptyset$ for each $ i\in k$ then
$F\cap \{q_i\}_{i\in k} \neq \emptyset$.
\end{itemize}
Notice that $\cal F$ is a closed set in ${\cal
P}(\Poset_{\alpha}^{m,n}\cap \Model)$. Next, let $\cal G$ be the set
of all $F\in \cal F$ such that $F$ is $\Poset_{\alpha}^{m,n}$ generic
over $ \Model$ and $p^*\in G$. To see that $\cal G$ is a dense
$G_\delta$ relative to $\cal F$, it suffices to show that if $D$ is a
dense subset of $\Poset_{\alpha}^{m,n}$ then $\{F\in {\cal F} : F\cap
D\neq
\emptyset\}$ is a dense open set in $\cal F$. In order to verify this,
let ${\cal V}$  be an open set in $\cal F$. It may be assumed that
$${\cal V} = \{F \in {\cal F} : p \in F  \AND \{q_i\}_{i\in k}\cap F =
\emptyset\}$$ for some $p\in \Poset_{\alpha}^{m,n}$ and $\{q_i\}_{i\in k}
\subseteq  \Poset_{\alpha}^{m,n}$. Because ${\cal V}$ is a non-empty
open set in $\cal F$ it must be that there is some $r \in
\Poset_{\alpha}^{m,n}$ such that
$r \leq p$ and $r \wedge q_i = \emptyset$ for each $ i\in k$. Let $r'
\in D$ be such that $r' \leq r$ and observe that that $\{F\in {\cal F}
: r'\in F\} \subseteq {\cal V}$. Hence ${\cal G}$ is a dense $G_\delta$
in a closed subspace of $\Poset_{\alpha}^{m,n}$.
It is worth noting that the third
condition in the definition of $\cal F$ implies that $\cal G$ has a
base consisting of sets of the form ${\cal U}(p) = \{F\in G : p \in F\}$. 

Let $\Psi :{\cal G}\mapsto \Reals^m$ be the mapping defined by
$\Model[G]\models x^G =
\Psi(G)$ where $x^G$ is the interpretation of the name $x$ in
$\Model[G]$.  Notice that $\Psi$ is continuous.

Before continuing, some notation concerning fusion arguments will be
established. Let $\{f_i\}_{i\in\omega}$ be a sequence of finite functions
satisfying the following properties:
\begin{itemize}
\item $f_i : \alpha \cap \Model \to \omega$
\item for each $i\in \omega$ there exists a unique ordinal $\beta(i)
\in \alpha \cap \Model$ such that $f_{i+1}(\gamma) = f_i(\gamma)$
unless $\gamma = \beta(i)$ and $f_{i+1}(\beta(i)) = f_i(\beta(i))+1$
where $f_i(\beta(i))$ is defined to be 0 if $\beta(i)$ is not in the
domain of $f_i$
\item $\sup_{i\in\omega}f_i(\gamma) = \omega$ for each $\gamma\in
\alpha\cap \Model$.
\end{itemize}
For each condition $p \in \Poset_{\alpha}^{m,n}$ the definition of
$\Sacks(m,n)$ guarantees that, for each $\gamma \in \alpha \cap
\Model$, there is a name $\{C^\gamma_\xi : \xi 
\in {}^{\stackrel{\omega}{\smile}}(n+1)\}$ such that
$p\restriction\gamma$ forces that $\{C^\gamma_\xi : \xi 
\in {}^{\stackrel{\omega}{\smile}}g(n+1)\}$ is a witness to the fact
that  $p(\gamma) \in \Sacks(m,n)$ and each $C^\gamma_\xi$ is an open
ball whose radius and centre are rational. It will also be assumed that
$ C^\gamma_\emptyset = \Reals^m$ for each $\gamma$.  

Suppose that  $r \in \Poset_{\alpha}^{m,n}$, $\Gamma\in
[\alpha\cap \Model]^{\aleph_0}$ and $A$ is a function with domain
$\Gamma$ such that $A(\gamma)$ is a $\Poset_{\gamma}^{m,n} $
name for a member of $\Sacks(m,n)$ for each $\gamma$ in $\Gamma$. Then
$r[A]$ denotes the condition defined by 
$$r[A](\gamma) = \left\{\begin{array}{ll}
r(\gamma) &\IF \gamma\notin \Gamma\\
r(\gamma)\cap A(\gamma) &\IF \gamma\in \Gamma\end{array}\right.$$  
For $i\in\omega$, $p \in \Poset_{\alpha}^{m,n}$ and $h\in
\prod_{\gamma\in\alpha\cap
\Model}{}^{f_i(\gamma)}(n+1) $ the notion of what it means for $p$ to
be determined is
the same as in \cite{ba.la.sack}. In particular, if $\Gamma 
\in [\alpha\cap \Model]^{<\aleph_0}$ and $f:\Gamma \to \omega$ then let
$$T_f = \bigcup_\beta\in\alpha \prod_{\gamma\in \Gamma\cap \beta}
{}^{f(\gamma )}(n+1)$$ and note that $T_f$ is a tree under inclusion.
 Define a condition $p \in
\Poset_{\alpha}^{m,n}$ to be $f$-determined 
with respect to 
$\{C^\gamma_\xi : \xi 
\in {}^{\stackrel{\omega}{\smile}}(n+1)\}$ 
if there is a function $A$
with domain $T_f$ such that if $$ h \in  \prod_{\gamma\in \Gamma\cap
(\beta + 1)}{}^{f(\gamma )}(n+1)$$ then, for each $\gamma\in \Gamma\cap
(\beta + 1)$, $A(h)(\gamma)$ is an $m$-ball with
rational centre and radius such that $p[A(h)\restriction \beta]
\forces{\Poset_{\beta}^{m,n}}{C^\beta_{h(\beta)} = A(h)(\beta)}$ and,
furthermore, $A(h\restriction \beta) = A(h) \restriction (\beta)$ for
every $\beta$. 

By a fusion sequence will be meant a sequence of conditions
$\{p_i\}_{i\in\omega}$ such that
\begin{itemize}
\item for each $i\in\omega$ and $\gamma$ in the support of  $f_i$ 
there is a sequence $$\{C^{\gamma,i}_\xi : \xi 
\in {}^{\stackrel{\omega}{\smile}}(n+1)\}$$ such that
$p_i\restriction\gamma$ forces that $\{C^{\gamma,i}_\xi : \xi
\in {}^{\stackrel{\omega}{\smile}}(n+1)\}$ is a witness to the fact
that  $p_i(\gamma) \in \Sacks(m,n)$
\item each $p_i$ is $f_i$-determined with respect to
$\{C^{\gamma,i}_\xi : \xi 
\in {}^{\stackrel{\omega}{\smile}}(n+1)\}$ by a function  $A_i$
defined on $T_i = T_{f_i}$
\item  $p_{i+1}\restriction\gamma$ forces that $C^{\gamma,i}_\xi 
= C^{\gamma,i+1}_\xi$ for each $i\in \omega$, $\gamma\in \alpha\cap
\Model $ and $\xi$ such that $\card{\xi} \leq f_i(\gamma)$
\end{itemize}
Notice that, since the sets $T_i$ are all disjoint,
there is no ambiguity in letting $A= \bigcup_{i\in\omega} A_i$ and
saying that $A$ witnesses that   
$p_i$ is $f_i$-determined with respect to
$$\{C^{\gamma,i}_\xi : \xi 
\in {}^{\stackrel{\omega}{\smile}}(n+1)\}$$ rather than 
that $\{A_i\}_{i\in\omega}$ witnesses this. This situation will be
referred to by 
saying that $A$ witnesses that $\{p_i\}$ is a fusion sequence.
If $G\in \cal G$ then define $H_{i}(G)$ to be the unique, maximal
member of $T_i$ such that $p_i[A(H_i(G))]\in G$ provided that such a
maximal member exists at all. Finally, if $h$ is a
maximal element of $ T_i$ then let ${\cal U}(h)$ be the open set in
$\cal G$ defined by 
${\cal U}(h) ={\cal U}(p_i[A(h)])$
and note that $i$ is the uniquely determined as the only integer such
that $h$ is a maximal member of $T_i$. Of course, ${\cal U}(h)$ and
$H_i(G)$ are only defined in the context of a given fusion sequence
but, since this will always be clear,
it will not be added to the notation.

Two facts are worth noting. First,  it is easy to verify that the fusion of
such a sequence is in $\Poset_{\alpha}^{m,n}$. Second, if $G$ and $G'$
are in $\cal G$ and $H_i(G)$ and $H_i(G')$ are defined and equal for
all $i\in \omega$ then $G = G'$.

Two cases will now be considered depending on whether or not $\alpha$
is a limit ordinal. 

\begin{lemma}
If $\alpha $ is limit \labelx{limit}ordinal then there is a 1-smooth curve in
$\Reals^m$, $L$,  and a condition $ r \leq p^*$ such that $q
\forces{\Poset_{\alpha}^{m,n}}{x\in L}$. 
\end{lemma}
\proof
For each $G\in \cal G$ let $\tau(G)$ be the set of all $V \in P_{1,m}$
such that for all $\epsilon > 0$, $q\in G$ and $\beta \in \alpha$ there is
$G'\in \cal G$ such that $G'\cap \Poset_{\beta}^{m,n} = G \cap
\Poset_{\beta}^{m,n}$, $q\in G'$ and $\rho(V,\Psi(G) -\Psi(G')) < \epsilon$.  
The first thing to observe is that if $G \in {\cal G}$ then $\tau(G)
\neq \emptyset$. To see this, first observe that
$\tau(G) = \bigcap_{\beta\in\alpha\cap \Model}\bigcap_{q\in
G}\overline{\tau(G,\beta,q)}$ where $$\tau(G,\beta,q) =
\{\langle \Psi(G) - \Psi(G')\rangle : G'\cap \Poset_{\beta}^{m,n} = G \cap
\Poset_{\beta}^{m,n}, q\in G' \AND \Psi(G) \neq \Psi(G')\}$$
and, since $P_{1,m}$ is compact, it suffices to show that
if $\beta \in \alpha \cap \Model$ and $q\in G$ then
$\tau(G,\beta,q)\neq \emptyset$.  Since
$p^*\forces{\Poset_{\alpha}^{m,n}}{x\notin V[G\cap \Poset_{\beta}^{m,n}]}$ for
every $\beta \in
\alpha$ it follows that there must be $G'$ containing both $p^*$ and $q$ 
which is $\Poset_{\alpha}^{m,n}/(G\cap \Poset_{\beta}^{m,n})$ generic
over $\Model[G\cap\Poset_{\beta}^{m,n}]$ and such that the
interpretation of $x$ in $\Model[G\cap \Poset_{\beta}^{m,n}*G']$ is
different from the interpretation of $x$ in $\Model [G]$.

Since $\cal G$ and $P_{1,m}$ are both Polish spaces and $\tau$ is a
Borel subset of ${\cal G}\times P_{1,m}$, it is possible to appeal to
the von Neumann Selection Theorem to find a Baire measurable function
$\Delta :{\cal G}\to P_{1,m}$ such that $\Delta(G) \in \tau(G)$ for
each $G\in \cal G$. 
Let $W\subseteq \cal G$ be a dense $G_\delta$ such that
$\Delta\restriction W$ is continuous. Let $W_n$ be dense open sets
such that $W = \bigcap_{i\in\omega}W_n$. 

The next step will be to construct a fusion sequence $\{p_i\}$ 
and $A$ which witnesses this so that the following conditions are
satisfied
\begin{itemize}
\item if $h$ is a maximal member of $T_i$ then ${\cal U}(h) \subseteq
W_i$
\item if $h_0$ and $h_1$ are distinct maximal members of $T_i$ then
the images
$\Psi({\cal U}(h_0))$ and $ \Psi({\cal U}(h_1))$ have disjoint
closures
\item if $G$  belongs to $W$, $H_i(G) = H_i(G')$ and
$H_{i+1}(G) \neq H_{i+1}(G')$ then
$\rho( \Psi(G) -\Psi(G'), \Delta(G)) < 1/i$
\end{itemize}
First, it will be shown that if such a fusion sequence can be
constructed, then the proof is complete. In particular, let $L$ be the
image under $\Psi$ of the set of all $G\in \cal G$ such that $H_i(G)$
is defined for all $i \in \omega$. Let $r$ be the fusion of the
sequence $\{p_i\}_{i\in \omega}$ and note that
$q\forces{\Poset_{\alpha}^{m,n}}{x\in L}$. 

To see that $L$ is 1-smooth, suppose that $x = \Psi(G)\in L$ and let
$\epsilon > 0$. It suffices to show that there is $\delta > 0 $ such
that if $x' \in L$ and $\|x - x'\| < \delta$ then
$\rho(\Delta(G),x-x') < \epsilon$. To this end, let $i >
1/\epsilon$ and let $\delta$ be so small that if $h$ and $h'$ are
distinct maximal members of $T_i$ then the distance between
$\Psi({\cal U}(h))$ and $\Psi({\cal U}(h'))$ is greater than
$\delta$. Therefore, if $\|x- x'\| < \delta$  it follows that
$H_i(G') = H_i(G)$ where $G'$ is such that $\Psi(G') = x'$. Let $j\in
\omega$ be the greatest integer such that $H_j(G) = H_j(G')$.  The
first property guarantees that both $G$ and $G'$ belong to $W$ and,
since $H_{j+1}(G)
\neq H_{j+1}(G')$ it follows that $\rho(\Delta(G),x-x') = \rho(
\Psi(G)-\Psi(G'), \Delta(G)) < 1/j < 1/i < \epsilon$.

The next claim will be used in showing that the desired fusion
sequence can be constructed.
\begin{claim}
Suppose that $\epsilon > 0$, $q\in
\Poset_{\alpha}^{m,n}$ and $\beta
\in \alpha$. Then there exist conditions $\{q_i\}_{i=0}^{ n}$
such that 
\begin{enumerate}
\item $q_i \leq q$ for $i\in n+1$
\item $q_i \restriction \beta = q_{i'}\restriction\beta$ for $i$ and
$i'$ in $ n+1$  
\item  if $\{G_i\}_{i=0}^{n} \subseteq W$ and  $q_i \in G_i$ for each
$i \leq n$
$$\rho(\Psi(G_i)-\Psi(G_j),\Delta(G_i)) < \epsilon$$ so long as 
 $i\neq j$
\item $q_i\forces{\Poset_{\alpha}^{m,n}}{x\in E_i}$ for closed balls
$E_i$ such that $E_i\cap E_{i'}= \emptyset$
if $i \neq i'$
\end{enumerate}
\end{claim}

Assuming the claim, suppose that $p_i$ and $A_i$ of the fusion
sequence have been constructed satisfying the induction
requirements. Let $\{h_s\}_{s\in k}$ enemurate all the maximal
elements of $T_i$ and let $\mu\in \alpha\cap \Model$ contain the domain of
$f_{i+1}$. Let $q_s^{0,j} = p_i[A(h_s)]$ for each $s\in k$.  Proceed
by induction on $y$ to construct $q_s^{y,j}$ such that
\begin{itemize}
\item $q_s^{y,j} \leq q_s^{y-1,j}$ for each $y\leq k$
\item if $\mu' \leq \mu$ and $h_s\restriction \mu' = h_{s'}\restriction \mu'$ 
then  $q_s^{y,j}\restriction \mu' = q_{s'}^{y,j'}\restriction \mu'$
\item ${\cal U}( q_y^{y+1,j})\subseteq W_i$ for $y\leq
k$ and $ j\in n+1$  
\item if $q_y^{y+1,j}\in G\in W$, $ q_y^{y+1,j'}\in G'\in W$ and
$j\neq j'$ then 
$\rho(\Psi(G) - \Psi(G'),\Delta(G)) < 1/i$
\item  for all $y\leq k$ and $ j\leq n$
$$q_y^{y+1,j}\restriction \beta(i)
\forces{\Poset_{\beta(i)}^{m,n}}{ C^{\beta(i),i}_{h_y(\beta(i))\wedge
j} = A^{j,h_y}}$$
\item $q_y^{y+1,j}\forces{\Poset_{\alpha}^{m,n}}{x\in E_y^j}$ where
$\{E_y^j\}_{j\in n+1}$ is a pairwise disjoint collection of closed balls 
\item if $y\leq s$ the $q_s^{y,j} = q_s^{y,j'}$ for all $j$ and $j'$
in $n+1$
\end{itemize}
 To see that this induction can be carried out, suppose that
$\{q_s^{y,j}\}_{s\in k}$ have been defined. From the last induction
hypothesis it follows that there is some $q''$ such that $q'' =
q_y^{y,j}$ for all $j\in n+1$. Since $W_i$ is dense open it is
possible to find $q' \leq q''$ such that ${\cal U}( q')
\subseteq W_i$. Extend $q'$ to $q$ such that 
$$q\restriction
\beta(i)\forces{\Poset_{\beta(i)}^{m,n}}{C^{\beta(i)}_{h_y(\beta(i))\wedge
j} = A^j}$$ for each $j\leq n$.
Use the claim to find conditions
$\{q_y^{y+1,j}\}_{j\in n+1}$ such that
\begin{itemize}
\item $q_y^{y+1,j} \leq q$ for $i\in n+1$
\item $q_y^{y+1,j} \restriction \mu = q_y^{y+1,j'}\restriction\mu$
for $j$  and $j'$ in $ n+1$
\item  if $\{G_j\}_{j\in n+1} \subseteq W$ and  $q_y^{y+1,j} \in G_j$ for each
$j \in n+1$
then $$\rho(\Psi(G_j)-\Psi(G_{j'}),\Delta(G_j)) < 1/(i+1)$$ so long as 
 $j\neq j'$
\item $q_y^{y+1,j}\forces{\Poset_{\alpha}^{m,n}}{x\in E_y^j}$ for closed balls
$E_y^j$ such that $E_y^j\cap E_y^{j'}= \emptyset$
if $j \neq j'$
\end{itemize}
If $s\neq y$ let $b(s)$ be the least member of the domain  of $f_i$
such that $h_s(b(s))\neq h_y(b(s))$ and define $q_s^{y+1,j}$ to be the
least upper bound  of $ q_y^{y+1,j}\restriction b(s)$ and $q_s^{y,j}$.
It is easily verified that all of the induction requirements are
satisfied.

Now, if $h$ is a maximal member of $T_{i+1}$ then let $j(h)\leq n$ and
$h^*$ be a maximal member of $T_i$ such that $h(\beta(i)) =
h^*(\beta(i))\wedge j(h)$ and, if $\gamma
\neq \beta(i)$ then $h(\gamma) = h^*(\gamma)$. Let $q_h =
q_{h^*}^{k,j(h)}$. 
Let
$$A_{i+1}(h)(\gamma) = \left\{\begin{array}{ll}
A_i(h^*)(\gamma) &\IF \gamma \neq \beta(i)\\
A^{j(h),h^*} &\IF \gamma = \beta(i)\end{array}\right.$$
Define $p_{i+1}$ to be the join of all the
conditions $q_h$ as $h$ ranges over all maximal members of $T_{i+1}$.
Note that if $y= y'$ then, by the induction hypothesis on $p_i$,
$\Psi({\cal U}(h_y))$ and $\Psi({\cal U}(h_{y'}))$ have disjoint
closures so it may be assumed that $E^j_y\cap E^{j'}_{y'}= \emptyset$.
Thus the three requirements of the desired fusion sequence are satisfied.

All that remains to be done is to prove the claim. To this end, let
$q$ , $\epsilon$ and $\beta$ be given. 
Using the fact that $W$ is
dense, let $G^*\in W$ be arbitrary such that $q\in G$. Using the
continuity of $\Delta$ on $W$, find $q'\in G$ such that
if $G'\in W$ and $q'\in G'$ then $\rho(\Delta(G),\Delta(G')) <
\epsilon/2$
 and let $p^{-1} = q'$. 
First notice that it suffices to construct by
induction on $i$ a sequence
 $\{(E_i,E^i,G_i,p_i,p^i) \}_{i\in n+1}$ such that for each $ i\in n+1$
\begin{enumerate}
\item $p_i \leq p^{i-1}$
\item $p^i \leq p^{i-1}$
\item $G_i\cap \Poset_{\beta}^{m,n} = G^*\cap \Poset_{\beta}^{m,n}$ 
\item $p^i \in G^*$ for each $i$
\item $p_i \in G_i$ for each $i$
\item $E^i$ and $E_i$ are disjoint closed subsets of $\Reals^m$
\item $E^{i+1}\cup E_{i+1}\subseteq E^i$
\item  $\Psi({\cal U}(p^i))$ is contained in the
interior of $E^i$
\item  $\Psi({\cal U}(p_i))$ is contained in the
interior of $E_i$
\item $\rho(x-x',\Delta(G^*)) < \epsilon/2$ if  $x \in E^i$ and $x'\in E_i$
\end{enumerate}
The reason this suffices is that, having done so, using conditions (3)
and (5) it is possible to find a single $p \in G^*\cap
\Poset_{\beta}^{m,n}$  extending each
$p_i$. Let $q_i$ be the greatest lower bound of both $p_i$ and $p$. It
follows that  $q_i \leq p_i \leq p^{i-1} \leq p^{-1} = q$. Moreover,
$q_i\restriction \beta = p$ for each $i\in n+1$. 
Also,  it follows from conditions (6) and (7) that
the sets $\{E_i\}_{i\in n+1}$ are pairwise disjoint closed sets and
$p_i\forces{}{x\in E_i}$ by condition (9). Finally, suppose that
$i\neq j$, $\{G_i,G_j\} \subseteq W$,  $q_i \in G_i$ and   $q_i \in G_i$.
Then  $\rho(\Delta(G_i),\Delta(G^*)) < \epsilon/2$ because $q_i  \leq q'$. 
Hence 
$$\rho(\Psi(G_i)-\Psi(G_j),\Delta(G_i)) < \epsilon$$ by condition (10).

To carry out the induction, suppose that
$(E_i,E^i,G_i,p_i,p^i) \}_{i\in J}$ have been constructed.
  From the
definition of $\Delta(G^*)$  it follows that there is $G_J$ such that
\begin{itemize}
\item $G_J\cap \Poset_{\beta}^{m,n} = G^*\cap \Poset_{\beta}^{m,n}$
\item  $p^{J-1}\in G_J$ 
\item $\rho(\Psi(G^*)-\Psi(G_J),\Delta(G^*)) < \epsilon/4$
\item $\Psi(G^*) \neq \Psi(G_J)$
\end{itemize}
Let $E^J$ and $E_J$ be disjoint closed neighbourhoods of $\Psi(G^*)$
and $\Psi(G_J)$ 
respectively such that $\rho(x-x',\Delta(G^*)) < \epsilon/2$ for
any $x \in E^J$ and $x'\in E_J$.  
 Since induction hypothesis (8) implies that
 $\Psi({\cal U}(p^{J-1}))$ is contained in the
interior of $E^{J-1}$ and
$ p^{J-1}\in G_J$ it follows that it may be assumed that
 $E^J\cup E_{J}\subseteq E^{J-1}$.
From the continuity of $\Psi$ it is
possible to find $p_J\in G_J$ and $p^J\in G^*$ extending $p^{J-1}$
  such that $\Psi({\cal U}(p_J))$ is contained in the interior of $ E_J$ and
$\Psi({\cal U}(p^J))$ is contained in the interior of $ E^J$.
All of the induction hypotheses are now satisfied
\stopproof

The possibility that $\alpha$ is a sucesssor must now be considered.
The proof has the same structure as the limit case but the details are
different. Only the successot case requires the use  of higher
dimensional tangent planes.

\begin{lemma}
If $\alpha $ is a successor \labelx{successor} ordinal then there is
an $n$-smooth surface in $\Reals^m$, $L$, and a condition $ r \leq p^*$
such that $q
\forces{\Poset_{\alpha}^{m,n}}{x\in L}$. 
\end{lemma}
\proof Let $\alpha = \beta + 1$. 
For each $G\in \cal G$ let $\tau(G)$ be the set of all $V \in P_{n,m}$
such that for all $\epsilon > 0$ and $q\in G$ there 
is a family $a \in [\cal G]^{n+1}$ such that 
\begin{enumerate}
\item $G\in a$
\item $a$ is $\epsilon$-orthogonal under some indexing
\item $G'\cap \Poset_{\beta}^{m,n} = G \cap
\Poset_{\beta}^{m,n}$ for each $G'\in a$
\item  $q\in \cap a$
\item $\Psi\restriction a$ is one-to-one
\item  $\rho(V,\Psi(G') - \Psi(G'')) < \epsilon$ for  $\{G',G''\}\in
[a]^2$
\end{enumerate}
  As in Lemma~\ref{limit}, it must be noted that if $G \in {\cal G}$
then $\tau(G)
\neq \emptyset$.  
Notice that $$\tau(G) \supseteq \bigcap_{q\in G}\bigcap_{\epsilon >
0}\overline{\tau(G,q,\epsilon)}$$ where $\tau(G,q,\epsilon)$ is
defined to be the set of all spaces generated by $$\{\Psi(G'') -
\Psi(G') : \{G',G''\}\in [a]^2\} $$ where $a$ satisfies condition (1)
to (5) with respect to $q$ and $\epsilon$.
Since $P_{n,m}$ is compact, it suffices to show that
if  $q\in G$ then
$\tau(G,q,\epsilon)\neq \emptyset$.  In $\Model[G\cap \Poset_{\beta}^{m,n}]$
the name $q(\beta) $ is interpreted as a condition in $\Sacks(m,n)$.
It follows that in $\Model[G\cap \Poset_{\beta}^{m,n}]$ 
 there must be open sets $\{C_i\}_{i=0}^{ n}$
such that $C_i\cap q(\beta) \in \Sacks(m,n)$ and
for any selection $x_i\in C_i$ the family $a= \{x_i\}_{i=0}^{ n}$ is
$\epsilon$-orthogonal. Since $$p^*\forces{\Poset_{\alpha}^{m,n}}{x\notin
\Model[G\cap \Poset_{\beta}^{m,n}]}$$ it is an easy matter to  
choose $G'_i$ which is $\Sacks(m,n)$ generic over
 $\Model[G\cap\Poset_{\beta}^{m,n}]$, $q_i \in G'_i$ and such that the
interpretation of $x$ in $\Model[G\cap \Poset_{\beta}^{m,n}*G'_i]$ is
different from the interpretation of $x$ in 
$\Model[G\cap \Poset_{\beta}^{m,n}*G'_j]$ unless $i=j$. It follows
that, letting $G_i = G\cap \Poset_{\beta}^{m,n}* G'_i$,
the space generated by  $$\{\Psi(G) - \Psi(G_i) \}_{i=1}^n  $$
belongs to $\tau(G,q,\epsilon)$.

As in Lemma~\ref{limit}, it is possible to find a Baire measurable function
$\Delta :{\cal G}\to P_{n,m}$ such that $\Delta(G) \in \tau(G)$ for
each $G\in \cal G$. 
Let $W\subseteq \cal G$ be a dense $G_\delta$ such that
$\Delta\restriction W$ is continuous. Let $W_n$ be dense open sets
such that $W = \bigcup_{i\in\omega}W_n$. 

The next step will be to construct a fusion sequence $\{p_i\}$ 
and $A$, which witnesses this, so that the following conditions are
satisfied
\begin{itemize}
\item if $h$ is a maximal member of $T_i$ then ${\cal U}(h) \subseteq
W_i$
\item if $h_0$ and $h_1$ are distinct maximal members of $T_i$ then
the images
$\Psi({\cal U}(h_0))$ and $ \Psi({\cal U}(h_1))$ have disjoint
closures
\item if $G$ belongs to $W$, $H_i(G) = H_i(G')$ and
$H_{i+1}(G) \neq H_{i+1}(G')$ then
$\rho(\Psi(G)-\Psi(G'), \Delta(G)) < 1/i$
\end{itemize}
First, it will be shown that if such a fusion sequence can be
constructed, then the proof is complete. In particular, let $L$ be the
image under $\Psi$ of the set of all $G\in \cal G$ such that $H_i(G)$
is defined for all $i \in \omega$. Let $r$ be the fusion of the
sequence $\{p_i\}_{i\in\omega}$. Note that $\Psi^{-1}\restriction L$
is defined.

The proof  that $L$ is $n$-smooth is the same as in Lemma~\ref{limit}. 
Suppose that $x = \Psi(G)\in L$ and let
$\epsilon > 0$. It suffices to show that there is $\delta > 0 $ such
that if $x' \in L$ and $\|x - x'\| < \delta$ then
$\rho(\Delta(G),x-x') < \epsilon$. To this end, let $i >
1/\epsilon$ and let $\delta$ be so small that if $h$ and $h'$ are
distinct maximal members of $T_i$ then the distance between
$\Psi({\cal U}(h))$ and $\Psi({\cal U}(h'))$ is greater than
$\delta$. Therefore, if $\|x- x'\| < \delta$  it follows that
$H_i(G') = H_i(G)$ where $G'$ is such that $\Psi(G') = x'$. Let $j\in
\omega$ be the greatest integer such that $H_j(G) = H_j(G')$.  The
first property guarantees that $G$  belongs to $W$ and,
since $H_{j+1}(G)
\neq H_{j+1}(G')$, it follows that $\rho(\Delta(G),x-x') = \rho(
 \Delta(G),\Psi(G) - \Psi(G')) < 1/j < 1/i < \epsilon$.

To see that the fusion sequence can be found, suppose that
$p_i$ and $A_i$ of the fusion sequence have been constructed
satisfying the induction requirements. Let $\{h_s\}_{s\in k}$
enumerate all the maximal elements of $T_i$. Let $q_s^{0,j} = p_i[A(h_s)]$ for
each $s\in k$.  Proceed by induction on $y$ to construct $q_s^{y,j}$
such that
\begin{itemize}
\item $q_s^{y,j} \leq q_s^{y-1,j}$ for each $y\leq k$
\item if $\mu \leq \alpha$ and $h_s\restriction \mu = h_{s'}\restriction \mu$ 
then  $q_s^{y,j}\restriction \mu = q_{s'}^{y,j'}\restriction \mu$
\item ${\cal U}( q_y^{y+1,j})\subseteq W_i$ for $y\leq
k$ and $ j\in n+1$  
\item if $q_y^{y+1,j}\in G\in W$, $ q_y^{y+1,j'}\in G'\in W$ and $j\neq j'$ then
$\rho(\Psi(G) - \Psi(G'),\Delta(G)) < 1/i$
\item  for all $y\leq k$ and $ j\leq n$
$q_y^{y+1,j}\restriction \beta(i)
\forces{\Poset_{\beta(i)}^{m,n}}{ C^{\beta(i),i}_{h_y(\beta(i))\wedge
j} = A^{j,h_y}}$ 
\item $q_y^{y+1,j}\forces{\Poset_{\alpha}^{m,n}}{x\in E_y^j}$ where
$\{E_y^j\}_{j\in n+1}$ is a pairwise disjoint collection of closed balls 
\item if $y\leq s$ the $q_s^{y,j} = q_s^{y,j'}$ for all $j$ and $j'$
in $n+1$
\end{itemize}
 To see that this induction can be carried out, suppose that
$\{q_s^{y,j}\}_{s\in k}$ have been defined. From the last induction
hypothesis it follows that there is some $q''$ such that $q'' =
q_y^{y,j}$ for all $j\in n+1$. Since $W_i$ is dense open it is
possible to find $q' \leq q''$ such that ${\cal U}( q')
\subseteq W_i$. Extend $q'$ to $q$ such that 
$$q\restriction
\beta(i)\forces{\Poset_{\beta(i)}^{m,n}}{C^{\beta(i)}_{h_y(\beta(i))\wedge
j} = A^{j,h_y}}$$ for each $j\leq n$. 

Using the fact that $ W$ is dense in $\cal G$, let $G\in W$ be such
that $q\in G$.  Using the definition of $\tau(G)$, let $\bar{\epsilon}
< 1/(i+1)$ and let $a\in [{\cal G}]^{n+1}$ be
$\bar{\epsilon}$-orthogonal such that
\begin{itemize}
\item $G\in a$
\item $G'\cap \Poset_{\beta}^{m,n} = G \cap
\Poset_{\beta}^{m,n}$ for each $G'\in a$
\item  $q\in \cap a$
\item $\Psi\restriction a$ is one-to-one
\item    $\rho(\Delta(G),\Psi(G') - \Psi(G'')) < \epsilon$ for  $\{G',G''\}\in
[a]^2$
\end{itemize}
It is therefore possible to choose $\hat{q}_i\in G_i$ as well as
$A_{i,h_y}$ such that
  $\hat{q}_i\restriction \beta = \hat{q}_0\restriction \beta$
and such that
$\hat{q}_0\forces{\Poset_{\beta}^{m,n}}{A_{i,h_y} = \hat{q}_i(\beta)}$
and, moreover, if $ z_i \in A_{i,h_y}$ 
for each $i \leq n$ then 
$\{z_i\}_{i=0}^{ n}$
is $frac{1}{i+1}$-orthogonal. Using the continuity of $\Delta$ and $\Psi$ as
well as the fact  that $p^*\forces{\Poset_{\alpha}^{m,n}}{x\notin
\Model[G\cap \Poset_{\beta}^{m,n}]}$ it is possible to extend
each $\hat{q}_j$ to $q_y^{y+1,j}$ such that
\begin{itemize}
\item $q_y^{y+1,j} \leq q$ for $i\in n+1$
\item $q_y^{y+1,j} \restriction \beta = q_y^{y+1,j'}\restriction\beta$
for $j$  and $j'$ in $ n+1$
\item  if $\{G_j\}_{j\in n+1} \subseteq W$ and  $q_y^{y+1,j} \in G_j$ for each
$j \in n+1$
then $$\rho(\Psi(G_j)-\Psi(G_{j'}),\Delta(G_j)) < 1/(i+1)$$ so long as 
 $j\neq j'$
\item $q_y^{y+1,j}\forces{\Poset_{\alpha}^{m,n}}{x\in E_y^j}$ for closed balls
$E_y^j$ such that $E_y^j\cap E_y^{j'}= \emptyset$
if $j \neq j'$
\end{itemize}
If $s\neq y$ let $b(s)$ be the least member of the domain  of $f_i$
such that $h_s(b(s))\neq h_y(b(s))$ and define $q_s^{y+1,j}$ to be the
least upper bound  of $ q_y^{y+1,j}\restriction b(s)$ and $q_s^{y,j}$.
It is easily verified that all of the induction requirements are
satisfied.

Now, if $h$ is a maximal member of $T_{i+1}$ then let $j(h)\in n+1$ and
$h^*$ be a maximal member of $T_i$ such that $h(\beta(i)) =
h^*(\beta(i))\wedge j(h)$ and, if $\gamma
\neq \beta(i)$ then $h(\gamma) = h^*(\gamma)$. Let $q_h =
q_{h^*}^{k,j(h)}$. 
Let
$$A_{i+1}(h)(\gamma) = \left\{\begin{array}{ll}
A_i(h^*)(\gamma) &\IF \gamma \neq \beta(i)\\
A^{j(h),h^*} &\IF \gamma = \beta(i) \neq \beta\\
{A}_{j(h),h^*} &\IF \gamma = \beta(i) = \beta
\end{array}\right.$$
Define $p_{i+1}$ to be the join of all the
conditions $q_h$ as $h$ ranges over all maximal members of $T_{i+1}$.
This satisfies the three requirements of the desired fusion sequence.

\stopproof
 Lemmas \ref{successor} and \ref{limit} together show that the ground
model smooth sets are sufficient to cover all real added by
iteratively adding reals with $\Sacks(m,n)$. Combined with
Corollary~\ref{GoUp}, this immediately gives the following theorem.

\begin{theor}
If  $1 \leq  n < m \in \omega$, then it is consistent, relative to the consistency
of set theory itself, that $\cov({\frak D}_{m,n}) < \cov({\frak
D}_{m,n+1})$.
\end{theor}

\section{Decomposing continuous functions}

In \cite{cimopaso} the authors consider the following question: If $\cal{A}$
and $\cal{B}$ are two families of functions between Polish spaces,
what is the least cardinal $\kappa$ such that every member of
$\cal{A}$ can be decomposed into $\kappa$ members of $\cal{B}$ --- this
cardinal $\kappa$ the authors call $\mbox{dec}(\cal{A},\cal{B})$.
The most natural class to consider for $\cal{B}$ is the class of
continuous functions and the problem to which it gives rise  had been
posed by Lusin who wondered whether every Borel function could be
decomposed in countably many continuous functions.
In more recent times, it has been shown by Abraham and Shelah that every
function of size less than $\frak c$ can be decomposed into countably
many continuous --- and even monotonic --- functions \cite{ab.ru.sh}.
 Various results concerning $\mbox{dec}(\cal{B}_1,\cal{C})$
where ${\cal B}_1$ is the class of pointwise limits of continuous
functiosn and $\cal C$ is the class of continuous functions can be
found in  \cite{step.34},\cite{step.30} and \cite{sole}.

The question of determining $\mbox{dec}(\cal{C},\cal{D})$ where $\cal
D$ is the class of functions which are differentiable on their domain
was raised by  M. Morayne and J. Cichon. It was not known whether it
is consistent that $\mbox{dec}(\cal{C},\cal{D}) < 2^{\aleph_0}$ and
the best lower bound for $\mbox{dec}(\cal{C},\cal{D})$ was noted by
Morayne to be the additivity of the null ideal. The following result
improves the lower bound and shows that 
$\mbox{dec}(\cal{C},\cal{D}) < 2^{\aleph_0}$ is indeed consistent.

\begin{theor} 
$\mbox{dec}(\cal{C},\cal{D}) = \cov({\cal D}_{2,1})$\labelx{Equiv}.
\end{theor}
\proof To show that
$\mbox{dec}(\cal{C},\cal{D}) \leq \cov({\cal D}_{2,1})$ let $\kappa = 
\mbox{dec}(\cal{C},\cal{D})$ and let $\{S_\eta\}_{\eta\in \kappa}$ de
a decomposition of $\Reals^2$ into 1-smooth curves. Given a continuous
function $f:\Reals \to \Reals$, let $f_\eta f\cap S_\eta$.

To see that $f_\eta$ is differentiable on its domain, let $x$ be in
the domain of $f_\eta$. Let $m$ be the slope of the tangent line to
$S_\eta$ at the point $(x,f_\eta(x))$ allowing the possibility that $m
= \infty$. Then, using the fact that $f$ is continuous, it is  
 easily verified that $f'_\eta(x) = m$.

To show that
$\mbox{dec}(\cal{C},\cal{D}) \geq \cov({\cal D}_{2,1})$ let $A
\subseteq \Reals$  be a perfect set and $\theta \in (0,\pi/4)$ be such that 
if $a$ abd $a'$ are distinct points in $A\times A$ then angle formed
by the horizontal axis and the line connecting $a$ and $a'$ is
different from $\theta$. Let $H_\theta$ be the function which projects
$\Reals^2$ to the vertical axis along the line at angle $\theta$ woth
respect to the horizontal axis. Let $H_{\pi/2}: \Reals\to \Reals$ be
the orthogonal projection onto the horizontal axis. Since
$H_\theta\restriction (A\times A)$ is one-to-one and $A\times A$ is
compact, it follows that $H_\theta^{-1}$ is continuous and, hence, so
is $H = H_{\pi/2}\circ H_\theta^{-1}$. Since the domain of $H$ is
comapct, it can be extended to a continuous function on th eentire
real line. 

Now suppose that $\{X_\eta\}_{\eta \in \kappa}$ are subsets of
$\Reals$ such that $H\restriction X_\eta$ is differentiable for each
$\eta \in \kappa$. To see that $H_\theta^{-1}(X_\eta)$ is a 1-smooth
curve, suppose that $(x,y)\in H_\theta^{-1}(X_\eta)$ is a point at
which $H_\theta^{-1}(X_\eta)$ does not have a 1-dimensional tangent.
In this case it is possible to find distinct slopes $m_0$ and $m_1$
and sequences $\{(x^i_n,y^i_n)\}_{n\in\omega}$ for each $i\in 2$ such
that $$\lim_{n\to\infty}\frac{y - y^i_n}{x- x^i_n} = m_i$$ allowing
the possibility of an infinite limit. Let $w^i_n = H^{-1}_\theta
(x^i_n,y^i_n)$ and $w = H_\theta^{-1}(x,y)$. It follows from the
linearity of $H_\theta$ that there are distinct $M_0$ and $M_1$ such
that $$\lim_{n\to\infty}\frac{x - x^i_n}{w - w^i_n} = M_i$$ and this
contradicts the differentiability of $H$ at $w$ since $H(w^i_n) =
x^i_n$ and $H(w) = x$.
\stopproof

\makeatletter \renewcommand{\@biblabel}[1]{\hfill#1.}\makeatother
\renewcommand{\bysame}{\leavevmode\hbox to3em{\hrulefill}\,}


\begin{thebibliography}{1}

\bibitem{ab.ru.sh}
A.~Abraham, M.~Rubin, and S.~Shelah, {\em On the consistency of some partition
  theorems for continuous colorings, and structure of $\aleph_1$-dense real
  order types}, Ann. Pure Appl. Logic {\bf 29} (1985), 123--206.

\bibitem{ba.la.sack}
J.~E. Baumgartner and R.~Laver, {\em Iterated perfect set forcing}, Ann. Math.
  Logic {\bf 17} (1979), 271--288.

\bibitem{cimopaso}
J.~Cicho\'{n}, M.~Morayne, J.~Pawlikowski, and S.~Solecki, {\em Decomposing
  {B}aire functions}, J. Symbolic Logic {\bf 56} (1991), no.~4, 1273--1283.

\bibitem{step.34}
S.~Shelah and J.~Stepr\={a}ns, {\em Decomposing baire class 1 functions into
  continuous functions}, Fund. Math. {\bf ??} (to appear), ?--?

\bibitem{sole}
S.~Solecki, Ph.D Thesis, 1994.

\bibitem{step.30}
J.~Stepr\={a}ns, {\em A very discontinuous {B}orel function}, J. Symbolic Logic
  {\bf 58} (1993), 1268--1283.

\end{thebibliography}
\end{document}